\documentclass{article}
\usepackage{amsfonts,amsbsy,graphics, graphicx}

\begin{document}
\pagestyle{empty}
{\sl{\small{\noindent Second International Workshop on\\ \noindent   Functional and Operatorial Statistics.\\ \noindent  Santander, June 16-18, 2011}}}\\~\\
\begin{center}{\Large{\bf Functional kernel estimators \\of conditional extreme quantiles}}
\end{center}
\begin{center}
\vspace{5mm}
{\bf{Laurent GARDES$^{**}$ and St\'ephane GIRARD$^{*}$}}\\
\vspace{5mm}
\noindent Team Mistis, INRIA Rh\^one-Alpes and LJK, Inovall\'ee, 655, av. de l'Europe,\\Montbonnot, 38334 Saint-Ismier cedex, France.\\
$^*$ {\tt Stephane.Girard@inrialpes.fr}

\vspace{5mm}
\rule{5cm}{.5mm}
\end{center}

\vspace{1cm}

\noindent {\bf\Large Abstract}\\

We address the estimation of ``extreme'' conditional quantiles {\it i.e.}
when their order converges to one as the sample size increases.
Conditions on the rate of convergence of their order to one
are provided to obtain asymptotically Gaussian distributed kernel estimators.
A Weissman-type estimator and kernel estimators of the conditional tail-index are derived, permitting
to estimate extreme conditional quantiles of arbitrary order.

\vspace{1cm}

\noindent {\bf\Large 1. Introduction}\\

Let $(X_i,Y_i)$, $i=1,\dots,n$ be independent
copies of a random pair $(X,Y)$ in $E \times{\mathbb{R}}$
where $E$ is a metric space associated to a distance $d$.
We address the problem of estimating 
$q(\alpha_n|x)\in{\mathbb{R}}$
verifying 
$
{\mathbb{P}}(Y>q(\alpha_n|x)|X=x)=\alpha_n
$
where $\alpha_n\to 0$ as $n\to\infty$ and $x\in E$.
In such a case, $q(\alpha_n|x)$ is referred to as an extreme
conditional quantile in contrast to classical conditional quantiles
(known as regression quantiles)
for which $\alpha_n=\alpha$ is fixed in $(0,1)$.
While the nonparametric estimation of ordinary regression quantiles 
has been extensively studied,
see for instance the seminal papers~(Roussas, 1969), (Stone, 1977) or (Ferraty and Vieu, 2006, Chapter~5)
less attention has been paid to extreme conditional quantiles despite their potential
interest. 
Here, we focus on the setting where the
conditional distribution of $Y$ given $X=x$ has an infinite endpoint
and is heavy-tailed,
an analytical characterization of this property being given in the next section.
We show, under mild conditions, that extreme conditional
quantiles $q(\alpha_n|x)$ can still be estimated through a functional kernel estimator
of ${\mathbb{P}}(Y>.|x)$. We provide sufficient conditions on the rate
of convergence of $\alpha_n$ to 0 so that our estimator is
asymptotically Gaussian distributed.
Making use of this, some functional kernel estimators of the conditional tail-index
are introduced and a Weissman type estimator (Weissman, 1978) is derived, permitting
to estimate extreme conditional quantiles $q(\beta_n|x)$ where $\beta_n\to 0$
arbitrarily fast.\\

\noindent {\bf\Large 2. Notations and assumptions}\\

The conditional survival function (csf) of $Y$ given $X=x$ is
denoted by $\bar{F}(y|x)={\mathbb{P}}(Y>y|X=x)$. 
The kernel estimator of $\bar{F}(y|x)$ is defined for all $(x,y)\in E \times{\mathbb{R}}$
by
\begin{equation}
\label{defestproba}
\hat{\bar{F_n}}(y|x)= \left.\sum_{i=1}^n K(d(x,X_i)/h) Q((Y_i-y)/\lambda) \right/ 
\sum_{i=1}^n K(d(x,X_i)/h) ,
\end{equation}
with $Q(t)=\int_{-\infty}^t q(s) ds$ where 
$K:{\mathbb{R}}^+\to{\mathbb{R}}^+$ and $q:{\mathbb{R}}\to{\mathbb{R}}^+$ are two kernel functions, and $h=h_n$ and $\lambda=\lambda_n$ are two nonrandom sequences such that 
$h\to 0$ as $n\to\infty$.
In this context, $h$ and $\lambda$ are called window-width.
This estimator was considered for instance in~(Ferraty and Vieu, 2006, page~56).
In Theorem~1, the asymptotic distribution of~(\ref{defestproba})
is established when estimating small tail probabilities, {\it i.e}
when $y=y_n$ goes to infinity with the sample size $n$.
Similarly, the kernel estimators of conditional quantiles 
$q(\alpha|x)$
are defined via the generalized inverse of $\hat{\bar{F_n}}(.|x)$: 
\begin{equation}
\label{defestquant}
\hat q_n(\alpha | x ) =\inf\{t,\; \hat{\bar{F_n}}(t|x)\leq \alpha\},
\end{equation}
for all $\alpha\in(0,1)$. 
Many authors are interested in this type of estimator for fixed $\alpha\in(0,1)$: weak and strong consistency are proved respectively in (Stone, 1977) and (Gannoun, 1990), asymptotic normality being established when $E$ is finite dimensional by (Stute, 1986), (Samanta, 1989), (Berlinet {\it et al.}, 2001)
and by (Ferraty {\it et al.}, 2005) when $E$ is a general metric space.
In Theorem~2,
the asymptotic distribution of~(\ref{defestquant})
is investigated when estimating extreme quantiles, {\it i.e}
when $\alpha=\alpha_n$ goes to 0 as the sample size $n$ goes to infinity.
The asymptotic behavior of such estimators depends on the nature of
the conditional distribution tail. In this paper, we focus on heavy tails.
More specifically, we assume that the csf satisfies\\

(A1): $\displaystyle\bar{F}(y|x)=c(x)  \exp\left\{-\int_1^y \left( \frac{1}{\gamma(x)}-\varepsilon(u|x)\right)\frac{du}{u}\right\}$,\\

\noindent where $\gamma(.)$ is a positive function of the covariate $x$,
$c(.)$ is a positive function and $|\varepsilon(.|x)|$ is continuous and ultimately 
decreasing to 0. (A1)~implies that the conditional distribution
of $Y$ given $X=x$ is in the Fr\'echet maximum domain of attraction.
In this context, $\gamma(x)$ is referred to as the conditional tail-index 
since it tunes the tail heaviness of the conditional distribution of
$Y$ given $X=x$. Assumption (A1)~also yields that
$\bar{F}(.|x)$ is regularly varying at infinity with index $-1/\gamma(x)$.
 {\it i.e} for all $\zeta>0$,
\begin{equation}
\label{defl}
 \lim_{y\to\infty} \frac{\bar{F}(\zeta y | x)}{\bar{F}(y|x)} = \zeta^{-1/\gamma(x)}.
\end{equation}
The function $\varepsilon(.|x)$ plays
an important role in extreme-value theory since it drives the
speed of convergence in~(\ref{defl}) and more generally the bias
of extreme-value estimators. Therefore, it may be of interest to
specify how it converges to~0. In (Gomes {\it et al.}, 2000),
the auxiliary function is supposed to be regularly varying
and the estimation of the corresponding regular variation index is addressed.
Some Lipschitz conditions are also required:\\

(A2): There exist $\kappa_\varepsilon$, $\kappa_c$, $\kappa_\gamma>0$ and $u_0>1$ such that for all $(x,x')\in E^2$ and $u>u_0$,
\begin{eqnarray*}
\displaystyle \left|\log c(x)-\log c(x')\right|&\leq &\kappa_c d(x,x'),\\
\displaystyle \left|\varepsilon(u|x)-\varepsilon(u|x')\right|&\leq& \kappa_\varepsilon d(x,x'),\\
\displaystyle \left|{1}/{\gamma(x)}-{1}/{\gamma(x')}\right|
&\leq &\kappa_\gamma d(x,x').
\end{eqnarray*}

\noindent The last assumptions are standard in the kernel estimation framework.\\

(A3): $K$ is a function with support $[0,1]$
and there exist $0<C_1<C_2<\infty$ such that $C_1\leq K(t) \leq C_2$ for all $t\in[0,1]$.

(A4): $q$ is a probability density function (pdf) with support $[-1,1]$.\\

\noindent One may also assume without loss of generality that $K$ integrates to one.
In this case, 
$K$ is called a type I kernel, see~(Ferraty and Vieu, 2006, Definition~4.1).
Finally, let $B(x,h)$ be the ball of center $x$ and radius $h$.
The small ball probability of $X$ is defined by $\varphi_x(h)={\mathbb{P}}(X\in B(x,h))$.
Under (A3), for all $\tau>0$, the $\tau$th moment is defined by $\mu^{(\tau)}_x(h)={\mathbb{E}}\{K^\tau(d(x,X)/h)\}$.\\

\noindent {\bf\Large 3. Main results}\\


Let us first focus on the estimation of small tail probabilities
$\bar{F}(y_n|x)$ when $y_n\to\infty$ as $n\to\infty$. 
The following result provides sufficient conditions for 
the asymptotic normality of $\hat{\bar{F_n}}(y_n|x)$.\\

\noindent {\bf{Theorem 1.}}
{\it Suppose (A1) -- (A4)~hold. Let $x\in E$ such that $\varphi_x(h)>0$ and  introduce $y_{n,j}=a_j y_n$ for $j=1,\dots, J$ with
$0<a_1<a_2< \dots < a_J$ and where $J$ is a positive integer.
If $y_n\to\infty$ such that  $n\varphi_x(h)  \bar{F}(y_n|x) \to \infty$,
$n\varphi_x(h)  \bar{F}(y_n|x)(\lambda/y_n)^2\to 0$ and $n \varphi_x(h) \bar{F}(y_n|x) (h \log y_n)^2  \to 0$ as $n\to \infty$, then
$$
\left\{ \sqrt{n  \mu_x^{(1)}(h)\bar{F}(y_n|x)}\left(\frac{\hat{\bar{F_n}}(y_{n,j}|x)
}{\bar{F}(y_{n,j}|x)}-1\right)\right\}_{j=1,\dots,J}
$$
is asymptotically Gaussian, centered, with covariance matrix 
$C(x)$ where $C_{j,j'}(x) = a^{1/\gamma(x)}_{j\wedge j'}$ for $(j,j')\in \{1,\dots, J\}^2$.\\
}

\noindent Note that $n \varphi_x(h) \bar{F}(y_n|x) \to \infty$
is a necessary and sufficient condition for 
the almost sure presence of at least one sample point
in the region $B(x,h)\times (y_n,\infty)$ of $E \times {\mathbb{R}}$.
Thus, this natural condition
states that one cannot estimate small tail probabilities out of the sample
using $\hat{\bar{F_n}}$. This result may be compared to (Einmahl, 1990) which 
establishes the asymptotic behavior of the empirical survival function
in the unconditional case but without assumption on the distribution.
Letting $\sigma_n(x)=(n \mu_x^{(1)}(h)\alpha_n)^{-1/2}$, the asymptotic normality
of $\hat q_n(\alpha_n|x)$ when $\alpha_n\to 0$ as $n\to \infty$ can be established under similar conditions.\\

\noindent {\bf{Theorem 2.}}
{\it
Suppose (A1) -- (A4)~hold. Let $x\in E$ such that $\varphi_x(h)>0$ and  introduce $\alpha_{n,j}=\tau_j \alpha_n$ for $j=1,\dots, J$ with
$\tau_1>\tau_2> \dots > \tau_J>0$ and where $J$ is a positive integer.
If $\alpha_n\to 0$ such that $\sigma_n(x) \to 0$, $ \sigma_n^{-1}(x)\lambda/q(\alpha_n|x)\to0$ and
$ \sigma_n^{-1}(x) h \log\alpha_n \to 0$ as $n\to \infty$,
then
$$
\left\{ \sigma_n^{-1}(x)\left(\frac{\hat q_n(\alpha_{n,j}|x)
}{q(\alpha_{n,j}|x)}-1\right)\right\}_{j=1,\dots,J}
$$
is asymptotically Gaussian, centered, with covariance matrix 
$
{\gamma^2(x)} \Sigma
$
where $\Sigma_{j,j'}= 1/\tau_{j\wedge j'}$ for $(j,j')\in \{1,\dots, J\}^2$.\\
}

\noindent The functional kernel estimator of extreme quantiles
$\hat q_n(\alpha_n|x)$ requires a stringent condition on the order $\alpha_n$
of the quantile, since by construction it cannot extrapolate beyond the maximum
observation in the ball $B(x,h)$. To overcome this limitation, a Weissman type
estimator (Weissman, 1978) can be derived:
$$
\hat q_n^{\mbox{\tiny W}}(\beta_n|x)= \hat q_n(\alpha_n|x) (\alpha_n/\beta_n)^{\hat\gamma_n(x)}.
$$
Here, $\hat q_n(\alpha_n|x)$ is the functional kernel estimator of the extreme quantile
considered so far and $\hat\gamma_n(x)$ is a functional estimator of the conditional
tail-index $\gamma(x)$. As illustrated in the next theorem, the extrapolation
factor $(\alpha_n/\beta_n)^{\hat\gamma_n(x)}$ allows to estimate extreme quantiles
of order $\beta_n$ arbitrary small.\\

\noindent {\bf Theorem~3.}
{\it
Suppose (A1)--(A4)~hold.
Let us introduce 
\begin{itemize}
\item $\alpha_n\to 0$ such that $\sigma_n(x) \to 0$, $ \sigma_n^{-1}(x)\lambda/y_n\to0$ and
$\sigma_n^{-1}(x)h\log\alpha_n \to 0$ as $n\to \infty$,
\item $(\beta_n)$ such that $\beta_n/\alpha_n\to 0$ as $n\to \infty$,
\item $\hat\gamma_n(x)$ such that $\sigma_n^{-1}(x)(\hat\gamma_n(x)-\gamma(x)) \stackrel{d}{\longrightarrow}{\mathcal {N}}(0,v^2(x))$ where $v^2(x)>0$.
\end{itemize}
Then, for all $x\in E$,
$$
\frac{\sigma_n^{-1}(x)}{\log(\alpha_n/\beta_n)}\left(\frac{\hat q_n^{\mbox{\tiny W}}(\beta_n|x)} {q(\beta_n|x)}-1\right)\stackrel{d}{\longrightarrow}{\mathcal{N}}(0,v^2(x)).
$$
}

\noindent Note that, when $K$ is the pdf of the uniform distribution,
this result is consistent with (Gardes {\it et al.}, 2010, Theorem~3),
obtained in a fixed-design setting.\\

\noindent Let us now give some examples of functional estimators of the conditional tail-index.
Let $\alpha_n\to 0$ and $\tau_1>\tau_2> \dots > \tau_J>0$ where $J$ is a positive integer.
Two additionary notations are introduced for the sake of simplicity:
$u=(1,\dots,1)^t \in {\mathbb{R}}^J$ and $v=(\log(1/\tau_1),\dots,\log(1/\tau_J))^t \in {\mathbb{R}}^J$.
The following family of estimators is proposed
$$
\hat\gamma_n(x) = \frac{\varphi (\log \hat q_n(\tau_1\alpha_n|x),\dots,\log \hat q_n(\tau_J\alpha_n|x))}
{\varphi(\log(1/\tau_1),\dots,\log(1/\tau_J))},
$$
where $\varphi:{\mathbb{R}}^J\to{\mathbb{R}}$ denotes a twice differentiable function
verifying the shift and location invariance conditions $\varphi( \theta v)=\theta \varphi(v)$
for all $\theta> 0$ and 
$\varphi( \eta u + x)= \varphi(x)$
for all $\eta\in{\mathbb{R}}$ and $x\in{\mathbb{R}}^J$.  
For instance, introducing the auxiliary function $m_p(x_1,\dots,x_J)=\sum_{j=1}^J (x_j-x_1)^p$ for all $p>0$ and considering 
$
\varphi_{\mbox{\tiny H}}(x)= m_1(x)
$
gives rise to a kernel version of the Hill estimator (Hill, 1975):
$$
{\hat\gamma_n^{\mbox{\tiny H}}(x)}= \sum_{j=1}^J \left[\log \hat q_n(\tau_j \alpha_n|x) -\log \hat q_n(\alpha_n|x)\right] \left / \sum_{j=1}^J \log (1/\tau_j)\right..
$$
Generalizations of the kernel Hill estimator can be obtained with $\varphi(x)=m_p(x)/m_1^{p-1}(x)$, see (Gomes and Martins, 2001, equation~(2.2)),
$\varphi(x)=m_p^{1/p}(x)$, see e.g. (Segers, 2001, example~(a)) or $\varphi(x)=m_{p\theta}^{1/\theta}(x)/m_{p-1}(x)$, $p\geq 1$, $\theta>0$,
see (Caeiro and Gomes, 2002).
In the case where $J=3$, $\tau_1=4$, $\tau_2=2$ and $\tau_3=1$, the function
$$
\varphi_{\mbox{\tiny P}}(x_1,x_2,x_3)= \log\left(\frac{\exp x_2 - \exp x_1}{\exp x_3 - \exp x_2}\right)
$$
leads us to a kernel version of Pickands estimator (Pickands, 1975)
$$
{\hat\gamma_n^{\mbox{\tiny P}}(x)}=\frac{1}{\log 2}\log\left(
\frac{\hat q_n(\alpha_n|x) - \hat q_n(2\alpha_n|x)}
{\hat q_n(2\alpha_n|x) - \hat q_n(4\alpha_n|x)}\right).
$$
We refer to~(Gijbels and Peng, 2000) for a different variant
of Pickands estimator in the context where the distribution of $Y$ given $X=x$ has a finite
endpoint. The asymptotic normality of $\hat\gamma_n(x)$ is a consequence of Theorem~2.\\

\noindent {\bf Theorem 4.}
{\it
Under assumptions of~Theorem~2 and if
$\sigma_n^{-1}(x)\varepsilon(q(\tau_1\alpha_n|x)|x)\to 0$ as $n\to \infty$, then,
$\sigma_n^{-1}(x) (\hat\gamma_n(x)-\gamma(x))$ converges to a centered Gaussian
random variable with variance
$$
V(x)=\frac{\gamma^2(x)}{\varphi^2(v)}(\nabla\varphi(\gamma(x)v))^t \Sigma  (\nabla\varphi(\gamma(x)v)).
$$
}

\noindent As an illustration, in the case of the kernel Hill and Pickands estimators, we obtain
\begin{eqnarray*}
V_{\mbox{\tiny H}}(x)&=& \gamma^2(x)\left(\sum_{j=1}^J \frac{2(J-j) +1}{\tau_j} - J^2\right) \left /
\left( \sum_{j=1}^J \log (1/\tau_j)\right)^2\right..\\
V_{\mbox{\tiny P}}(x)&=& \frac{\gamma^2(x) (2^{2\gamma(x)+1}+1)}{4(\log 2)^2(2^{\gamma(x)}-1)^2}.
\end{eqnarray*}
Clearly, $V_{\mbox{\tiny P}}(x)$ is
the variance of the classical Pickands estimator, see for instance (de Haan and Ferreira, 2006, Theorem~3.3.5). 
Focusing on the kernel Hill estimator and 
choosing $\tau_j=1/j$ for each $j=1,\dots,J$ yields $V_{\mbox{\tiny H}}(x)=\gamma^2(x)J(J-1)(2J-1)/(6\log^2(J!))$.
In this case, $V_{\mbox{\tiny H}}(x)$ is a convex function of $J$ and is minimum for $J=9$ leading to $V_{\mbox{\tiny H}}(x)\simeq 1.25\gamma^2(x)$.

\noindent {\bf\Large References}\\

\noindent A. Berlinet, A. Gannoun and E. Matzner-L\o ber.
Asymptotic normality of convergent estimates of conditional quantiles.
\textit{ Statistics}, 35:139--169, 2001.

\noindent
F. Caeiro and M.I. Gomes.
Bias reduction in the estimation of parameters of rare events. {\em Theory of Stochastic Processes}, 8:67--76, 2002.

\noindent
J.H.J Einmahl.
The empirical distribution function as a tail estimator,
{\it Statistica Neerlandica}, 44:79--82, 1990.

\noindent
F.~Ferraty and P.~Vieu.
\newblock {\em Nonparametric functional data analysis}.
\newblock Springer, 2006.

\noindent
F.~Ferraty, A. Rabhi and P.~Vieu.
Conditional quantiles for dependent functional data with application to the climatic {\it El Nino} Phenomenon,
{\it Sankhya: The Indian Journal of Statistics}, 67(2):378--398, 2005.

\noindent
A. Gannoun. Estimation non param\'etrique de la m\'ediane conditionnelle, m\'e\-dianogramme et m\'ethode du noyau,  \textit{Publications de l'Institut de Statistique de l'Universit\'e de Paris},  XXXXVI:11--22, 1990.

\noindent
L.~Gardes, S. Girard and A. Lekina.
\newblock Functional nonparametric estimation of conditional extreme quantiles.
\newblock {\em Journal of Multivariate Analysis}, 101:419--433, 2010.

\noindent
I.~Gijbels and L.~Peng.
\newblock {Estimation of a support curve via order statistics.}
\newblock {\em Extremes},
 3:251--277, 2000.

\noindent
M.I. Gomes and M.J. Martins and M. Neves.
Semi-parametric estimation of the second order parameter, asymptotic and finite
sample behaviour. {\em Extremes}, 3:207--229, 2000.

\noindent
M.I. Gomes and M.J. Martins.
Generalizations of the Hill estimator -  asymptotic versus finite
sample behaviour. {\em Journal of Statistical Planning and Inference}, 93:161--180, 2001.

\noindent
L. de Haan and A. Ferreira.
{\em Extreme Value Theory: An Introduction},
Springer Series in Operations Research and Financial Engineering,
Springer, 2006.

\noindent
B.M. Hill.
A simple general approach to inference about the tail of a distribution.
{\em The Annals of Statistics}, 3:1163--1174, 1975.

\noindent
J.~Pickands. Statistical inference using extreme
order statistics. {\it The Annals of Statistics}, 3:119--131, 1975.

\noindent
G.G.~Roussas. Nonparametric estimation of the transition
distribution function of a Markov process.
{\it Ann. Math. Statist.}, 40:1386--1400, 1969.

\noindent
T.~Samanta. Non-parametric estimation of conditional
 quantiles. \textit{Statistics and Probability Letters}, 7:407--412, 1989.

\noindent
J.~Segers. Residual estimators.
{\it Journal of Statistical Planning and Inference}, 98:15--27, 2001.

\noindent
C.J. Stone.  Consistent nonparametric regression (with discussion). \textit{ The Annals of Statistics}, 5:595--645, 1977.

\noindent
W.~Stute. Conditional empirical processes. \textit{The Annals of Statistics}, 14:638--647, 1986.

\noindent
I. Weissman. Estimation of parameters and large quantiles based on the $k$ largest observations, {\it{Journal of the American Statistical Association}}, 73:812--815, 1978.
\end{document}